\documentclass[12pt, reqno]{amsart}
\usepackage{amssymb,latexsym, amsmath, amscd, array}

\numberwithin{equation}{section} 
\numberwithin{figure}{section}
\numberwithin{table}{section}

\DeclareMathOperator{\syscat}{{\rm cat_{sys}}}
\DeclareMathOperator{\cat}{{\mbox{\rm cat$_{\rm LS}$}}}

\DeclareMathOperator{\vol}{{\rm vol}}

\DeclareMathOperator{\stsys}{{\rm stsys}}
\DeclareMathOperator{\sys}{{\rm sys}} 

\DeclareMathOperator{\sysh}{{\rm sysh}}

\DeclareMathOperator{\pisys}{{\rm sys}\pi}

\def\rp{\mathbb R\mathbb P}

\def\la{\langle}
\def\ra{\rangle} 
\def\ov{\overline} 
\def\ds{\displaystyle}

\newcommand{\N}{\mathbb N}

\newcommand{\R}{\mathbb R} 
\newcommand{\RP}{{\mathbb R}{\mathbb P}}
\newcommand{\T}{\mathbb T}
\newcommand{\Z}{\mathbb Z}
\newcommand{\gmetric}{{\mathcal G}}
\newcommand{\hmetric}{{\mathcal H}}

\newtheorem{theorem}{Theorem}[section]

\newtheorem{proposition}[theorem]{Proposition}
\newtheorem{lemma}[theorem]{Lemma}
\newtheorem{corollary}[theorem]{Corollary}

\newtheorem{cor}[theorem]{Corollary}

\theoremstyle{definition} 
\newtheorem{definition}[theorem]{Definition}
\newtheorem{example}[theorem]{Example}
\newtheorem{remark}[theorem]{Remark}

\newtheorem{question}[theorem]{Question}

\newcommand\lemref{Lemma~\ref}
\newcommand\propref{Proposition~\ref}
\newcommand\corref{Corollary~\ref}

\newcommand{\regular}{{B}}

\def\ie {{\it i.e.\ }}  \def\eg {{\it e.g.\ }}
\def\cf {\hbox{\it cf.\ }}

\def\dist{{\rm dist}}

\def\ga{\alpha} 
 
\def\gl{\lambda} 

\def\eps{\varepsilon} 
\def\gf{\varphi}

\def\m{\medskip}

\date{11 june 2007}
\begin{document}

\author[M.~Katz]{Mikhail G. Katz$^{*}$} \address{Department of
Mathematics, Bar Ilan University, Ramat Gan 52900 Israel}
\email{katzmik@math.biu.ac.il} \thanks{$^{*}$Supported by the Israel
Science Foundation (grants no.~84/03 and 1294/06)}

\author[Y.~Rudyak]{Yuli B. Rudyak$^{**}$} \address{Department of
Mathematics, University of Florida, PO Box 118105, Gainesville, FL
32611-8105 USA} \email{rudyak@math.ufl.edu} \thanks{$^{**}$Supported
by NSF, grant 0406311}

\title[Bounding volume by systoles of~$3$-manifolds] {Bounding volume
by systoles of~$3$-manifolds}



\begin{abstract}
We prove a new systolic volume lower bound for non-orientable
$n$-manifolds, involving the stable 1-systole and the codimension 1
systole with coefficients in~$\Z_2$.  As an application, we prove that
Lusternik-Schnirelmann category and systolic category agree for
non-orientable closed manifolds of dimension 3, extending our earlier
result in the orientable case.  Finally, we prove the homotopy
invariance of systolic category.
\end{abstract}

\maketitle

\tableofcontents

\section{Introduction}

The systolic project in its modern form was initiated by
M. Gromov~\cite{Gr1}, when he proved a volume lower bound for a closed
essential Riemannian manifold~$M$, which is curvature-free and depends
only on the least length of a noncontractible loop in~$M$, \ie the
1-systole~$\pisys_1(M)$:
\begin{equation}
\label{11z}
(\pisys_1(M))^n \leq C_n \vol_n(M).
\end{equation}
See \cite{bcg2} for a recent application.  Alternative approaches to
the proof of \eqref{11z} may be found in \cite {AK, We, Gu}.

More generally, one considers higher dimensional systoles, and seeks
analogous volume lower bounds.  The defining text for this material is
the monograph~\cite{Gr3} (which is an extended English version of
\cite{Gr0}), with additional details in the earlier texts
\cite{Gr1,Gr2}.  Recently there has been a considerable amount of
activity related to systolic inequalities.  The Loewner inequality and
its generalisations are studied in \cite{Am, IK, KL, BCIK1, BCIK2,
KS1, KRS}.  Near-optimal asymptotic bounds are studied in \cite{Ba3,
bal, BB, KS2, KSV}.  See \cite{SGT} for an overview of systolic
problems.

\section{Motivation}

The notion of systolic category~$\syscat$ was introduced by the
authors in \cite{KR1}, \cf \eqref{21b} below.  It can be regarded as a
differential-geometric analogue of the Lusternik-Schnirelmann category
(LS category)~$\cat$, \cf \cite{LS}.  While the definitions of the two
categories use completely different language, the values of the two
invariants turn out to be very close in a number of interesting cases.
Thus, the insight gained from the study of one of them can be used to
glean a fresh perspective on the other.

The significance of systolic category as compared to LS category can
be illustrated by the recent progress on a long-standing conjecture
concerning the latter \cite{DKR}.  Such progress was motivated in part
by systolic geometry.  Namely, in collaboration with A.~Dranishnikov,
we proved the 1992 conjecture of Gomez-Larra\~naga and
Gonzalez-Acu\~na, to the effect that~$n$-manifolds ($n\geq 3$) of LS
category~$2$ necessarily have free fundamental group.

Conversely, the knowledge of the value of LS category
for~$4$-manifolds
$M$ with non-free fundamental groups lends plausibility to the
possible existence of systolic inequalities on~$M$ corresponding to
the partition~$4=1+1+2$, \cf \eqref{2.1}.  Such an inequality is
currently known to exist only in a limited number of cases, due
essentially to Gromov \cite[Theorem~7.5.B]{Gr1}.

We will use the modern definition of the Lusternik-Schnirelmann
category, \cf \cite{CLOT}, which differs by a unit from the original
definition.  Thus,~$\cat$ of a contractible space is equal to zero.
The two invariants share a number of characteristics, including lower
bound by cup-length and sensitivity to Massey products, \cf
\cite[Theorem~11.1]{KR1}, \cite{Ka4}.  

\begin{remark}
The lowest known dimension of a manifold whose systolic category is
strictly smaller than its LS category, is dimension~$16$, see
\cite[Example~9.5]{KR1} or \cite[p.~104-105]{SGT}.
\end{remark}

We prove that systolic category and LS category agree for
non-orientable closed manifolds of dimension 3, extending our earlier
result in the orientable case \cite[Corollary~6.2]{KR1}.  The required
lower bound for the systolic category of a non-orientable manifold
follows from a new inequality~\eqref{12} involving systoles of
dimension and codimension one.  A different but similar inequality in
the orientable case was studied in~\cite{BK1, BK2}.  The proof
exploits harmonic forms, the Cauchy-Schwartz inequality, and the
coarea formula, see Section~\ref{proof}.

The paper is organized as follows.  In Section~\ref{one}, we recall
the definition of the systolic invariants.  In Section~\ref{thms}, we
recall the definition of the Berg\'e-Martinet constant, and present an
optimal systolic inequality~\eqref{12} combining systoles of dimension
and codimension 1, valid for non-orientable manifolds.  Its proof
appears in Section~\ref{proof}.  The positivity of systoles is proved
in Section~\ref{pos}.  In Section~\ref{five}, we prove the homotopy
invariance of systolic category, which parallels that of the LS
category.  Some open questions are posed in Section~\ref{new}.

All manifolds are assumed to be closed, connected, and smooth. All
polyhedra are assumed to be compact and connected, unless explicitly
mentioned otherwise.  To the extent that our paper aims to address
both a topological and a geometric audience, we attempt to give some
indication of proof of pertinent results that may be more familiar to
one audience than the other.

\section{Systoles and systolic category}
\label{one}

Let~$X$ be a (finite) polyhedron equipped with a piecewise Riemannian
metric~$\gmetric$.  We will now define the systolic invariants of
$(X,\gmetric)$.  Note that we adopt the usual convention that the
infimum calculated over an empty set is infinite.  Thus by definition,
a~$k$-systole is infinite in the case when~$H_k(X)=0$.

\begin{definition}
The homotopy 1-systole, denoted~$\pisys_1(X,\gmetric)$, is the least
length of a non-contractible loop in~$X$.  The homology 1-systole,
denoted~$\sysh_1(X,\gmetric)$, is defined in a similar way, in terms
of loops which are not zero-homologous.
\end{definition}

Clearly, we have~$\pisys_1\le \sysh_1$.

\begin{definition}
\label{32}
Let~$k\in \N$.  Higher
homology~$k$-systoles~$\sysh_k=\sysh_k(X,\gmetric, A)$, with
coefficients over a ring~$A=\Z$ or~$\Z_2$, are defined similarly
to~$\sysh_1$, as the infimum of~$k$-areas of~$k$-cycles, with
coefficients in~$A$, which are not zero-homologous.
\end{definition}

More generally, let~$\regular$ be the group of deck transformations of
a regular covering space of~$X$.  Then the homology groups of the
covering space of~$X$ can be identified with the homology groups
of~$X$ with coefficients in the group ring~$A[\regular]$.  Allowing
more general coefficients in such a group ring, we can therefore
define the corresponding systole
\begin{equation}
\label{ag}
\sysh_k(X, \gmetric; A[\regular]).
\end{equation}
Note that we adopt the usual convention, convenient for our purposes,
that the infimum over an empty set is infinity.

More detailed definitions appear in the survey \cite{CK} by C. Croke
and the first author, and in \cite{SGT}.  We do not consider higher
``homotopy'' systoles.

\begin{definition}[\cf \cite{Fe1, BK1}]
Given a class~$\ga\in H_k(X;\Z)$ of infinite order, we define the
stable norm~$\| \ga_\R \|$ by setting
$$
\|\ga_\R\|=\lim_{m\to \infty} m^{-1} \inf_{\ga(m)} \vol_k(\ga(m)),
$$
where~$\alpha_\R$ denotes the image of~$\alpha$ in real homology,
while~$\ga(m)$ runs over all Lipschitz cycles with integral
coefficients representing~$m\ga$. The stable homology~$k$-systole,
denoted~$\stsys_k(\gmetric)$, is defined by minimizing the stable norm
$\| \alpha_\R \|$ over all integral~$k$-homology classes~$\alpha$ of
infinite order.
\end{definition}

We have $\stsys_k \le \sysh_k$ in the absence of torsion.  The stable
systoles can be significantly different from the ordinary ones.  The
first examples of such a phenomenon, already for the~$1$-dimensional
homology systoles, were discovered by Gromov and described by
M. Berger in \cite{Be5}.  On the other hand, for~$\RP^n$ we
have~$\sysh_1 < \infty$ while~$\stsys_1=\infty$.

Recall that, in our convention, the systolic invariants are infinite
when defined over an empty set of loops or cycles.

\begin{remark}
M. Berger~\cite{B} defined invariants which eventually came to be
known as the~$k$-systoles, in the framework of Riemannian manifolds.
All systolic notions can be defined similarly for polyhedra, \cf
\cite{Gr2} and \cite{Ba02, Ba06}.  Note that every smooth manifold is
triangulable and therefore can be viewed as a polyhedron.  When~$k=n$
is the dimension,~$\sysh_n(M,\gmetric)$ is equal to the
volume~$\vol_n(M,\gmetric)$ of a compact
Riemannian~$n$-manifold~$(M,\gmetric)$.  For an~$n$-polyhedron~$X$,
however, the volume may not agree with the~$n$-systole~$\sys_n(X)$, as
the former is always finite, while the latter may be infinite,
when~$X$ does not possess a fundamental class. Moreover, it can happen
that~$\sys_nX\ne \vol_nX$ even if~$\sys_nX$ is finite: for example,
if~$X$ is a wedge of two~$n$-spheres.
\end{remark}

The idea is to bound the total volume from below in terms of
lower-dimensional systolic invariants.  Here we wish to incorporate
all possible curvature-free systolic inequalities, stable or unstable.
More specifically, we proceed as follows.

\begin{definition}
\label{21b}
\label{sys} Given~$k\in \N, k>1$ we set
\[
\ \sys_{k}(X, \gmetric)= \inf\{\sysh_k(X, \gmetric;\Z[\regular]),
\sysh_k(X, \gmetric;\Z_2[\regular]),\stsys_k(X, \gmetric) \},
\]
where the infimum is over all groups~$B$ of regular covering spaces of
$X$.  Furthermore, we define
$$
\sys_{1}(X, \gmetric)=\min \{\pisys_1(X, \gmetric),\stsys_1(X,
\gmetric)\}.
$$
\end{definition}

Note that the systolic invariants thus defined are positive (or
infinite), see Section~\ref{pos}.

Let~$X$ be an~$n$-dimensional polyhedron, and let~$d\geq 2$ be an
integer.  Consider a partition
\begin{equation}
\label{2.1} 
n= k_1 + \ldots + k_d,
\end{equation}
where~$k_i\geq 1$ for all~$i=1,\ldots, d$.  We will consider
scale-invariant inequalities ``of length~$d$'' of the following type:
\begin{equation}
\label{dd} \sys_{k_1}(\gmetric) \sys_{k_2}(\gmetric) \ldots
\sys_{k_d}(\gmetric) \leq C(X) \vol_n(\gmetric),
\end{equation}
satisfied by all metrics~$\gmetric$ on~$X$, where the constant~$C(X)$
is expected to depend only on the topological type of~$X$, but not on
the metric~$\gmetric$.  Here the quantity~$\sys_k$ denotes the infimum
of all non-vanishing systolic invariants in dimension~$k$, as defined
above.

\begin{definition}
\label{syscat} 
Systolic category of~$X$, denoted~$\syscat(X)$, is the largest
integer~$d$ such that there exists a partition~\eqref{2.1} with
\[
\prod\limits^{d}_{i=1} \sys_{k_i}(X,\gmetric) \leq C(X) \vol_n(X,\gmetric)
\]
for all metrics~$\gmetric$ on~$X$.  If no such partition and
inequality exist, we define systolic category to be zero.
\end{definition}

In particular,~$\syscat X \le \dim X$.  

\begin{remark}
Clearly, systolic category equals~$1$ if and only if the polyhedron
possesses an~$n$-dimensional homology class, but the volume cannot be
bounded from below by products of systoles of positive codimension.
Systolic category vanishes if~$X$ is contractible.  Another example of
a 2-polyhedron~$X$ with~$\syscat X=0$ is a wedge of the disk and the
circle, \cf \corref{zerosys}.
\end{remark}

\section{Inequality combining dimension and codimension 1}
\label{thms}

Given a maximal rank lattice~$L$ in a normed space~$(\R^b, ||\cdot
||)$, let~$\gl_1(L)$ denote the least length of a non-zero vector
of~$L$.

\begin{definition}
\label{11}
The Berg\'e-Martinet constant \cite{BM}, denoted~$\gamma_b'$, is
defined as follows:
\begin{equation}
\label{bm}
\gamma'_b = \sup\left\{ \lambda_1(L) \lambda_1(L^*)\left| L \subseteq
\R ^b \right. \right\},
\end{equation}
where the supremum is extended over all lattices~$L$ in~$\R^b$ with
its Euclidean norm.
\end{definition}

Here the dual lattice~$L^*\subset \R^b$ by definition consists of
elements~$y\in \R^b$ satisfying~$\la x,y \ra \in\Z$ for all~$x\in L$,
where~$\la\;,\;\ra$ is the Euclidean inner product.  A lattice
attaining the supremum in \eqref{bm} is called {\em dual-critical}.

Both the Hermite constant~$\gamma_b$ and the Berg\'e-Martinet constant
$\gamma'_b$ are asymptotically linear in~$b$.  Their values are known
in dimensions up to~$4$ as well as certain higher dimensions.

\begin{example}
\label{bmexample}
In dimension~$3$, the value of the Berg\'e-Martinet constant,~$\gamma
_3'= \sqrt{\frac{3}{2}}=1.2247\ldots$, is slightly below the Hermite
constant~$\gamma_3= 2^{\frac{1}{3}}=1.2599\ldots$.  It is attained by
the face-centered cubic lattice, which is not isodual
\cite[p.~31]{MH}, \cite[Proposition 2.13(iii)]{BM}, \cite{CS}.
\end{example}

We generalize an inequality proved in the orientable case in \cite
{BK1, BK2}.

\begin{theorem}
\label{33}
Let~$M$ be an~$n$-dimensional manifold with first Betti number~$b\geq
1$.  Then every metric~$\gmetric$ on~$M$ satisfies the systolic
inequality
\begin{equation}
\label{12}
\stsys_1(\gmetric) \sys_{n-1}(\gmetric;\Z_2) \leq \gamma_b'
\vol_n(\gmetric),
\end{equation}
where~$\gamma_b'$ is the Berg\'e-Martinet constant of \eqref{bm}.
Furthermore, inequality \eqref{12} is optimal.
\end{theorem}

This inequality is proved in Section~\ref{proof}.

\begin{corollary}\label{firstbetti}
We have~$\syscat(M) \ge 2$ for all manifolds~$M$ with positive first
Betti number.
\qed
\end{corollary}

\begin{corollary}\label{cor:cat=2}
If~$M$ is a closed~$3$-dimensional manifold with nontrivial free fundamental
group, then~$\syscat M=2$.
\end{corollary}

\begin{proof} 
The classifying space of the free group is a wedge of circles.  Thus a
manifold~$M$ with free fundamental group is not essential.  Therefore
its systolic category is at most~$2$.

The systolic inequality provided by \corref{firstbetti}, corresponding
to the partition~$3=1+2$ as in \eqref{2.1}, shows that~$\syscat M=2$.
\end{proof}

\begin{corollary}
\label{26}
Systolic category and LS category coincide for all closed
connected~$3$-manifolds, orientable or not.  \qed
\end{corollary}

\begin{proof}
The fact that the LS category of a~$3$-manifold is determined by its
fundamental group was proved by J.~Gomez-Larra\~naga and
F.~Gonzalez-Acu\~na \cite{GG} (see also~\cite{OR}). In particular,
$\cat M^3=3$ if and only if~$\pi_1(M)$ is not free.

By Gromov's inequality for essential manifolds \cite{Gr1}, combined
with I.~Babenko's converse to it \cite{Bab1}, an~$n$-manifold has
systolic category~$n$ if and only if it is essential.  Here the
systolic inequality corresponds to the partition~$n=1+1+\ldots+1$ as
in \eqref{2.1}.  In \cite[Corollary 7.3]{KR1}, we proved that every
$3$-manifold with non-free fundamental group is essential. Hence, if
$\cat (M^3)=3$ then~$M^3$ is essential, snd hence~$\syscat
M=3$. Furthermore, if~$\cat M=2$ then~$\pi_1(M)$ is nontrivial and
free, and hence~$\syscat M=2$ by \corref{cor:cat=2}. Finally, if~$\cat
M^3=1$ then~$M$ is a homotopy sphere, and thus~$\syscat M=1$.
\end{proof}

\begin{remark}
The class of 3-dimensional Poincar\'e complexes is essentially larger
than the class of 3-manifolds.  For example, the Sphere Theorem does
not hold for 3-dimensional Poincar\'e complexes by J.~Hillman's
work~\cite{Hi2}.  Hillman's example~$Y$ is irreducible, essential, and
virtually free.  Thus~$Y$ does not easily fit into the algebraic
dichotomy in the context of 3-manifolds discussed in \cite
[Proposition~7.2] {KR1}, \cf the ``Tits alternative'' of \cite{Hi1}.
It remains to be seen how the existence of such an example affects the
calculation of the two categories.
\end{remark}

\begin{question}
Does the conclusion of Corollary~\ref{26} hold more generally for
3-dimensional Poincar\'e complexes?  
\end{question}

\section{Proof of optimal inequality}
\label{proof}

With the proof of Theorem~\ref{33} in mind, let~$H_1(M; \Z)_\R$ be the
integer lattice in~$H_1(M;\R)$, and similarly for cohomology.  Given a
metric~$\gmetric$, one defines the stable norm~$\|\cdot\|$ in homology
and the comass norm~$\|\cdot\|^*$ in cohomology.  The normed lattices
\[
\left( H^{\phantom{1}}_1(M;\Z)_\R, \|\cdot\| \right) \mbox{ and } \left(
H^1(M;\Z)_\R, \|\cdot\|^* \right)
\]
are dual, whether or not~$M$ is orientable \cite[item~5.8]{Fe2}.  An
intuitive explanation for such duality may be found in
\cite[Proposition 4.35, p.~261]{Gr3}.

Let~$\|\cdot\|_2^*$ be the~$L^2$-norm in~$H^1(M;\R)$, \ie
\[
\|\omega\|_2^*=\inf_{\xi\in\omega}|\xi|_2
\]
where the infimum is over all closed forms~$\xi\in \omega$, and
$|\cdot|_2$ is the~$L^2$-norm for forms.  We have
$\|\omega\|_2^*=|\eta|_2$ for the harmonic form~$\eta \in \omega$
\cite{LM}, and in particular the norm~$\|\cdot\|_2^*$ is Euclidean.
We will consider the invariant~$\lambda_1(L) \lambda_1(L^*)$ for the
lattice~$L =H_1(M; \Z)_\R \subset H_1(M;\R)$.

\begin{lemma}\label{optimal}
Let~$\omega\in H^1(M;\Z)_\R$ be a cohomology class whose modulo~$2$
reduction~$\ov \omega \in H^1(M;\Z_2)$ is nonzero. Then
\begin{equation*}
\label{34}
\sys_{n-1}(\gmetric; \Z_2) \leq \|\omega\|^*_2
(\vol_n(\gmetric))^{1/2}.
\end{equation*}
\end{lemma}

\begin{proof}
Let~$\eta\in\omega$ be the harmonic 1-form for the metric~$\gmetric$
on~$M$.  Then~$\eta$ can be represented as~$df$ for some map
\[
f: M \to S^1=\R/\Z.
\]
Using the Cauchy-Schwartz inequality, we obtain
\begin{equation}
\label{31}
\begin{aligned}
\|\omega\|^*_2 (\vol_n(\gmetric))^{1/2} =
|\eta|_{2}(\vol_n(\gmetric))^{1/2}\geq \int_ {M } |df| d\vol_n ,
\end{aligned}
\end{equation}
where~$|\cdot|$ is the pointwise norm defined by the Riemannian metric.
We now use the coarea formula, \cf \cite[3.2.11]{Fe1},
\cite[p.~267]{Ch2}:
\begin{equation}
\label{24}
\int_ {M } |df| d\vol_n = \int_{S^1} \vol_{n-1} \left (
f^{-1}_{\phantom{I}}(t) \right)dt .
\end{equation}
Note that, for every regular value~$t$ of~$f$, the~$\Z_2$-homology
class of the hypersurface~$f^{-1}(t) \subset M$ is Poincar\'{e}
$\Z_2$-dual to~$\ov \omega$.  Hence,
\[
\vol_{n-1} \left (
f^{-1}_{\phantom{I}}(t) \right)\ge \sys_{n-1}(\gmetric; \Z_2)
\] 
for all regular values~$t$ of~$f$. By Sard's Theorem, the set of
regular values of~$f$ has measure~$1$ in~$S^1$. Thus,
\begin{equation}
\label{82}
\begin{aligned}
\int_ {M } |df| d\vol_n \geq \sys_{n-1}(\gmetric; \Z_2).
\end{aligned}
\end{equation}
The lemma results by combining inequalities~\eqref{31}
and \eqref{82}.
\end{proof}

\begin{proof}[Proof of Theorem~$\ref{33}$]
Let~$\|\cdot\|_2$ be the norm in homology dual to the~$L^2$-norm
$\|\cdot\|_2^*$ in cohomology.  Let~$\alpha\in H_1(M;\Z)_\R$ be an
element of least norm, so that
\[
\| \alpha \|_2 = \lambda_1 \left( H^{\phantom{1}}_1(M; \Z)_\R,
\|\cdot\|_2 \right) .
\]
 Clearly, 
$\ds\|\cdot\|_2^* \leq \|\cdot \|^* \vol_n(\gmetric)^{1/2}$, and so, dually, 
$\ds\|\alpha\| \leq \|\alpha\|_2 \vol_n(\gmetric)^{1/2}$. Choose~$\omega$ so 
that
\[
\|\omega \|^*_2 = \lambda_1(H^1(M;\Z)_{\R}, \|\cdot\|^*_2) .
\]
By \lemref{optimal}, we obtain
\[
\begin{aligned}
\stsys_1(\gmetric) \sys_{n-1}(\gmetric;\Z_2) & = \|\alpha \| \;
\sys_{n-1} (\gmetric;\Z_2) \\ & \leq \|\alpha \| \; \| \omega \|^*_2
\vol_n(\gmetric)^{1/2} \\ & \leq \|\alpha \|_2 \; \| \omega \|^*_2
\vol_n(\gmetric).
\end{aligned}
\]
The theorem now follows from the inequality
\[
\begin{aligned}
\|\alpha \|_2 \; \| \omega \|^*_2 &= \lambda_1 \left( H_1(M; \Z)_\R,
 \|\cdot \|_2^{\phantom{I}} \right) \lambda_1 \left( H^1(M; \Z)_\R, \|
 \cdot \|_2 ^* \right) \\&\leq \gamma_b'
\end{aligned}
\] 
by Definition~\ref{11} of the Berg\'e-Martinet constant.  The
optimality of the inequality results by considering a suitable product
metric on the product~$\T^b\times \RP^2$, where~$\T^b$ is dual-critical.
\end{proof}

\section{Positivity of systoles}
\label{pos}

\begin{proposition}\label{posit}
The homotopy~$1$-systole and the stable systoles defined in
Section~$\ref{one}$ are nonzero for all polyhedra~$X$.
\end{proposition}

\begin{proof}
We cover~$X$ by a finite number of open, contractible sets~$U_i$.  By
the Lebesgue Lemma, there exists~$\delta>0$ such that every subset of
$X$ of diameter at most~$\delta$ is contained in some~$U_i$.
Therefore the diameter of any non-contractible loop~$L$ must be more
than~$\delta$, and thus the length of~$L$ must exceed 2$\delta$,
proving the positivity of the 1-systole.

For the stable~$k$-systoles, the positivity follows by a
``calibration'' argument.  Namely, suppose classes~$\alpha\in
H_k(X;\Z)_{\R}$ and~$\omega\in H^k(X;\Z)_{\R}$ pair non-trivially and
positively.  (For a theory of differential forms on polyhedra see \eg
\cite{Ba02}.) Then
$$
1 \leq \int _\alpha \omega \leq C \| \omega \|^* \vol_k
(\alpha),
$$
and, moreover, 
$$
1 \leq \int _\alpha \omega \leq m^{-1}C \| \omega\|^* \vol_k(
\alpha (m))
$$  
for all~$m\in \N$. Minimizing over all singular Lipschitz cycles
$\alpha(m) \in m[\alpha_\R]$, we obtain the necessary bound~$1\leq C
\| \omega\|^*\| \alpha_{\R} \|~$. Now the result follows because the
abelian group~$H^k(X;\Z)_{\R}$ is of finite rank.
\end{proof}

For the ordinary~$k$-systoles, one cannot use differential forms, as
in the proof of \propref{posit}, due to possible torsion classes in
homology.  Nevertheless, the the positivity of~$\sysh_k$ holds as
well.  In fact, H.~Federer~\cite[item~4.2.2(1)]{Fe1} proved that
cycles with small mass are homologically trivial.  However, we need a
slightly stronger conclusion, to obtain a uniform lower bound for
systoles of covering spaces.  The following lemma is a consequence of
the construction used in the proof of the deformation theorem of
Federer and Fleming~\cite{FF, Fe1}, \cf \cite{Wh}.  The proof was
summarized in~\cite[Prop. 3.1.A]{Gr1}.

\begin{lemma}
\label{dist}
Let~$V$ be a~$k$-dimensional polyhedron in~$\R^N$.  Then there exists
a continuous map~$f$ of~$V$ into a~$(k-1)$-dimensional polyhedron
$K^{k-1}$ in~$\R^n$ such that
\[
\dist(v,f(v))\le C_N(\vol V)^{1/k}
\]
for all~$v\in V$ and for some constant~$C_N$ depending only on the
ambient dimension, where~$\dist$ denotes the Euclidean distance in
$\R^N$.  \qed
\end{lemma}

\begin{proposition}\label{pos2}
The homology~$k$-systoles~$\sysh_k(X,\gmetric;A)$ defined in
Section~$\ref{one}$, Definition~$\ref{32}$, are nonzero for all
polyhedra~$X$.  In fact, a uniform lower bound for~$\sysh_k$ is valid
for all covering spaces of~$X$.
\end{proposition}

\begin{proof}
We may view~$X$ as a polyhedron in~$\R^N$, since the metric~$\gmetric$
is bilipschitz equivalent to the restricted metric.  Note that~$X$ has
a regular neighborhood~$U\subset \R^N$, \ie~$X$ is a deformation
retract of~$U$.  Since~$X$ is compact, there exists~$\eps>0$ such that
$\dist(X,\R^N\setminus U)>\eps$. Let~$j: X \to U$ denote the
inclusion.
 
Note that we can regard any singular chain in~$X$, with coefficients
in~$A$ with~$A=\Z$ or~$\Z_p$, as a singular, not necessarily connected
polyhedron.  Consider a such a polyhedron~$\gf: F^k \to X$
representing a non-zero~$k$-dimensional homology with coefficients in
$A$.  We now view~$\gf$ as a singular polyhedron in~$\R^N$.

Choose~$\sigma>0$ such that~$C_N\sigma^{1/k}<\eps$ where~$C_N$ is the
constant from \lemref{dist}.  Assume that~$\vol(\gf(F))<\sigma$. Then
by Lemma~\ref{dist}, there is a map~$f: \gf(F) \to K^{k-1}\subset
\R^N$ such that~$\dist (a, f(a))<\eps$ for all~$a\in \gf(F)$.  Thus
$f(\gf(F))\subset U$.  Moreover, for all~$a\in \gf(F)$, the segment
joining~$a$ and~$f(a)$ is contained in~$U$. Therefore the maps
$$
\CD F @>\gf >> X @>j>>U\quad \text{ and }\quad F @>\gf >>\gf(F) @>f >>
U \endCD
$$
are homotopic.  But~$f(\gf(F))\subset K^{k-1}$, and so~$f\circ \gf : F
\to U$ represents the zero element in~$H_k(U;A)$. Hence~$j\circ \gf$
represents the zero element in~$H_k(U;A)$. But then~$\gf$ represents
the zero element in~$H_k(X;A)$. This is a contradiction.

The deformation of~$F$ into a polyhedron of positive codimension is
contained in~$U$.  Since~$X$ is a retract of~$U$, and because of the
Cellular Approximation Theorem, we can map the deformation into~$X$ in
a way that the last-moment map sends~$F$ into the~$(k-1)$-skeleton
$X^{(k-1)}$ of~$X$.

It follows that the same~$\delta$ works for the systoles of arbitrary
covering spaces of~$X$, by the covering homotopy property, \cf
\cite[Section~2] {BCIK2}.
\end{proof}

\section{Homotopy invariance of systolic category}
\label{five}

The homotopy invariance of systolic category follows from the
techniques developed by I. Babenko in \cite{Bab1}, \cf the compression
theorem of \cite[Theorem~5.1]{KR1}.  More precisely, we have the
following theorem, proved in \cite{Bab1} for the 1-systole by
essentially the same argument.

\begin{theorem}
The optimal systolic ratio associated with a partition of~$n=\dim(M)$
is a homotopy invariant of a closed manifold~$M$.
\end{theorem}

\begin{proof}
Let~$f: M^n \to N^n$ be a homotopy equivalence of closed PL manifolds.
By A. Wright \cite[Theorem 7.3]{Wr},~$f$ is homotopic to a PL monotone
map.  Recall that a continuous map is called monotone if the inverse
image of every point is connected and compact.  Thus, we can assume
the every~$n$-dimensional simplex of~$N$ has exactly one inverse image
simplex.  Then we can pull back systolic inequalities from one
manifold to the other, in the following sense.  The pullback metric
has {\em the same volume\/} as the target metric.  Meanwhile, the
projection map is distance decreasing, and therefore the target
manifold has smaller systoles than the source manifold.

Now suppose that~$M$ satisfies a systolic inequality 
\begin{equation}
\label{61b}
\prod_i \sys_{k_i}(M)\le C\vol_n(M),
\end{equation}
relative to a suitable partition~$n=k_1+\ldots+k_d$, for {\em all\/}
metrics, with a constant~$C$.  In particular, it satisfies it for all
pullback metrics.  Then the target manifold will satisfy the same
inequality (with the same partition of~$n$) with exactly the same
constant~$C$.  Thus the associated optimal systolic ratio, which is
the least such constant~$C$, is a homotopy invariant.
\end{proof}

\m Since Wright's result is not available for polyhedra, we have to
make do with a fixed simplicial map.  The disadvantage here is that
the pullback metric may have greater volume than the target, but
anyway it is controlled by the number of simplices in the inverse
image of a top dimensional simplex.  Thus the optimal systolic ratio
is no longer a homotopy invariant, but nevertheless the homotopy
invariance of systolic category persists, in the following sense.

\begin{theorem}
\label{monot}
Given two~$n$-dimensional polyhedra~$X$ and~$Y$, assume that there
exists a simplicial map~$f: X \to Y$ that induces an isomorphism in
$\pi_1$ and a monomorphism in homology with coefficients
in~$\Z[\regular]$ and~$\Z_2[\regular]$, where~$\regular$ runs over all
groups of regular covering maps of~$X$~$($and therefore~$Y)$. Then
$\syscat X \le \syscat Y$.
\end{theorem}

\begin{proof} 
The case~$\syscat X=0$ is trivial, so assume that~$\syscat
X>0$. Consider a partition~$n=k_1+\ldots+k_d$.  Suppose that~$X$
satisfies a systolic inequality 
\[
\prod_i \sys_{k_i}(X, \gmetric) \le C(X) \vol_n(X, \gmetric)
\]
for {\em all\/} piecewise Riemannian metrics.  Choose a piecewise
Riemannian metric on~$Y$ and consider the pull back (degenerate)
metric~$\hmetric$ induced by~$f$ on~$X$.  By the monomorphism
hypothesis, we have~$\sys_k(X,\hmetric) \ge\sys_k(Y)$. Meanwhile,
$\vol_n(X,\hmetric)\le k\vol_n(Y)$ where~$k$ is number of
$n$-simplices in~$X$.  A small perturbation of~$\hmetric$ will yield a
non-degenerate metric~$\gmetric$, satisfying the inequality
$\vol_n(X,\gmetric)\le 2k\vol_n(Y)$.  Moreover, we can also assume
that~$\sys_k(X,\gmetric)\ge\sys_k(Y)-\delta$ where~$\delta\le
\sys_k(X)$ for all~$k$.  Hence,
$$
\prod \sys_{k_i}(Y)\le \prod (\sys_{k_i}(X,\gmetric)+\delta)\le 
2C\vol_n(X,\gmetric)\le 4kC\vol_n(Y),
$$
and thus~$\syscat X \le \syscat Y$.
\end{proof}

We have the following two immediate corollaries.

\begin{corollary}\label{homottinv}
The systolic category of~$n$-dimensional polyhedra is a homotopy
invariant.  \qed
\end{corollary}

\begin{cor}\label{zerosys} 
Let~$X$ be an~$n$-dimensional polyhedron which is homotopy equivalent
to a polyhedron~$Y$ of dimension at most~$n-1$.  Then~$\syscat X=0$.
\end{cor}


\section{New directions}
\label{new}

\begin{question}
In the context of higher Massey products, is there a generalization of
the result \cite[Theorem~11.1]{KR1}, to a case where triple products
vanish but there is a higher nontrivial product?
\end{question}

\begin{question}
Considering finite covers which are intermediate between the free
abelian cover and the manifold.  Is it true that if the fiber class is
non-zero in the free abelian cover, then it is already nontrivial in a
finite cover?  In such case we would immediately get a lower bound of
$b_1(X)+1$ for~$\cat$, which parallels the bound for systolic category
resulting from \cite{IK}.
\end{question}

\begin{question}
Consider an (absolute) degree 1 map~$f : X \to T^2$ from a
non-orientable surface to the torus.  Consider also the product map
\[
        f : X \times S^2 \to T^2 \times S^2.
\]
The range has systolic category 3, by real cup length argument.
However, the domain has real cup length only 2.  Therefore we do not
have an immediate lower bound of 3 for systolic category, unlike LS
category.  Is there such a bound?  What is the systolic category of
$\rp^2 \times S^2$?
\end{question}

\begin{question}
Given a function~$f$ on a manifold, how does~$\syscat (f^c)$ of the
sublevel set~$f^c$ change as a function of~$c$?
\end{question} 

\section{Acknowledgment}

We are grateful to Brian White for a discussion of the deformation
theorem, exploited in Section~\ref{pos}.



\begin{thebibliography}{ABCDE}

\bibitem[AK00]{AK} Ambrosio, L.; Kirchheim, B.: Currents in metric
spaces.  {\em Acta Math.} \textbf{185} (2000), no.~1, 1--80.


\bibitem[Am04]{Am} Ammann, B.: Dirac eigenvalue estimates on two-tori.
{\em J. Geom. Phys.} \textbf{51} (2004), no. 3, 372--386.

\bibitem[Ba93]{Bab1} Babenko, I.: Asymptotic invariants of smooth
manifolds.  {\em Russian Acad. Sci. Izv. Math.} \textbf{41} (1993),
1--38.


\bibitem[Ba02]{Ba02} Babenko, I.: Forte souplesse intersystolique de
vari\'{e}t\'{e}s ferm\'{e}es et de poly\`{e}dres.  \textit{Annales de
l'Institut Fourier\/} \textbf{52} (2002), no. 4, 1259-1284.

\bibitem[Ba04]{Ba3} Babenko, I.: G\'eom\'etrie systolique des
vari\'et\'es de groupe fondamental~$\Z_2$, {\em
S\'emin. Th\'eor. Spectr. G\'eom. Grenoble}, \textbf{22} (2004),
25-52.



\bibitem[Ba06]{Ba06} Babenko, I.: Topologie des systoles
unidimensionelles. {L'Enseignement Math\'ematique (2)} \textbf{52}
(2006), 109-142.




\bibitem[BaB05]{BB} Babenko, I.; Balacheff, F.: G\'eom\'etrie
systolique des sommes connexes et des rev\^etements cycliques,
{\em Mathematische Annalen\/} (to appear).

\bibitem[Bal04]{bal} Balacheff, F.: Sur des probl\`emes de la
g\'eom\'etrie systolique. {\em S\'emin.  Th\'eor.  Spectr.  G\'eom.
Grenoble\/} \textbf{22} (2004), 71--82.




\bibitem[BaCIK05]{BCIK1} Bangert, V; Croke, C.; Ivanov, S.; Katz, M.:
Filling area conjecture and ovalless real hyperelliptic surfaces, {\em
Geometric and Functional Analysis (GAFA)\/} \textbf{15} (2005) no.~3,
577-597.  See \texttt{arXiv:math.DG/0405583}


\bibitem[BaCIK07]{BCIK2} Bangert, V; Croke, C.; Ivanov, S.; Katz, M.:
Boundary case of equality in optimal Loewner-type inequalities.  {\em
Trans. Amer. Math. Soc.} \textbf{359} (2007), no.~1, 1--17.  See
\texttt{arXiv:math.DG/0406008}




\bibitem[BaK03]{BK1} Bangert, V.; Katz, M.: Stable systolic
inequalities and cohomology products, {\em Comm. Pure Appl. Math.}
\textbf{56} (2003), 979--997.  Available at the site
\texttt{arXiv:math.DG/0204181}


\bibitem[BaK04]{BK2} Bangert, V; Katz, M.: An optimal Loewner-type
systolic inequality and harmonic one-forms of constant norm.  {\em
Comm. Anal. Geom.}  \textbf{12} (2004), number 3, 701--730.  See
\texttt{arXiv:math.DG/0304494}


\bibitem[BeM89]{BM} Berg\'{e}, A.-M.; Martinet, J.: Sur un probl\`eme
de dualit\'e li\'e aux sph\`eres en g\'eom\'etrie des nombres.  {\em
J.~Number Theory\/} \textbf{32} (1989), 14--42.

\bibitem[Be72]{B} Berger, M.: Du c\^ot\'e de chez Pu. {\em
Ann. Sci. \'Ecole Norm. Sup.} (1972), \textbf{4}, 1--44; A l'ombre de
Loewner, {\em ibid}, 241--265.


\bibitem[Be93]{Be5} Berger, M.: Systoles et applications selon Gromov.
S\'eminaire N.~Bourbaki, expos\'{e} 771, \textit {Ast\'{e}risque\/}
\textbf{216} (1993), 279--310.


\bibitem[BeCG03]{bcg2} Besson, G.; Courtois, G.;, Gallot, S.:
Hyperbolic manifolds, amalgamated products and critical exponents.
{\em C. R. Math. Acad. Sci. Paris\/} \textbf{336} (2003), no. 3,
257--261.


\bibitem[Ch93]{Ch2} Chavel, I.: Riemannian geometry -- a modern
introduction.  {\em Cambridge Tracts in Mathematics}, \textbf{108}.
Cambridge University Press, Cambridge, 1993.



\bibitem[CoS94]{CS} Conway, J. H.; Sloane, N. J. A.: On lattices
equivalent to their duals.  {\em J. Number Theory\/} \textbf{48}
(1994), no. 3, 373--382.

\bibitem[CLOT03]{CLOT} Cornea, O.; Lupton, G.; Oprea, J.; Tanr\'e, D.:
Lusternik-Schnirelmann category.  {\em Mathematical Surveys and
Monographs}, \textbf{103}.  American Mathematical Society, Providence,
RI, 2003.

\bibitem[CK03]{CK} Croke, C.; Katz, M.: Universal volume bounds in
Riemannian mani\-folds, {\it Surveys in Differential Geometry\/}
\textbf{VIII} (2003), 109--137.  Available at the site
\texttt{arXiv:math.DG/0302248}


\bibitem[DKR07]{DKR} Dranishnikov, A.; Katz, M.; Rudyak, Y.: Small
values of Lusternik-Schnirelmann and systolic categories for
manifolds. Available at the site arXiv:0706.1625


\bibitem[Fe69]{Fe1} Federer, H.: Geometric measure theory.  {\em
Grundlehren der mathematischen Wissenschaften},
\textbf{153}. Springer--Verlag, Berlin, 1969.

\bibitem[Fe74]{Fe2} Federer, H.: Real flat chains, cochains, and
variational problems.  \textit{Indiana Univ.\ Math.\ J.}  \textbf{24}
(1974), 351--407.

\bibitem[FF60]{FF} Federer, H.; Fleming, W.
Normal and integral currents. \textit{Ann. of Math.} (2) \textbf{72} (1960) 
458--520.

\bibitem[GG92]{GG} G\'omez-Larra\~naga, J.; Gonz\'alez-Acu\~na, F.:
Lusternik-Schnirel'mann category of~$3$-manifolds.  {\em Topology\/}
\textbf{31} (1992), no. 4, 791--800.


\bibitem[Gr81]{Gr0} Gromov, M.: Structures m\'etriques pour les
vari\'et\'es riemanniennes.  Edited by J. Lafontaine and
P. Pansu. {\em Textes Math\'ematiques}, \textbf{1}. CEDIC, Paris,
1981.

\bibitem[Gr83]{Gr1} Gromov, M.: Filling Riemannian manifolds, {\em
J. Diff. Geom.} \textbf{18} (1983), 1--147.

\bibitem[Gr96]{Gr2} Gromov, M.: Systoles and intersystolic
inequalities, Actes de la Table Ronde de G\'{e}om\'{e}trie
Diff\'{e}rentielle (Luminy, 1992), 291--362, {\em S\'{e}min. Congr.},
\textbf{1}, Soc. Math. France, Paris, 1996.

\bibitem[Gr99]{Gr3} Gromov, M.: Metric structures for Riemannian and
non-Riemannian spaces, {\em Progr. in Mathematics}, \textbf{152},
Birkh\"{a}user, Boston, 1999.


\bibitem[GrL83]{GL} Gromov, M.; Lawson, H. B., Jr.: Positive scalar
curvature and the Dirac operator on complete Riemannian manifolds.
{\em Inst. Hautes Etudes Sci. Publ. Math.}, \textbf{58} (1983),
83--196 (1984).


\bibitem[Gu06]{Gu}Guth, L.: Volumes of balls in large Riemannian
manifolds.  Available at \texttt{arXiv:math.DG/0610212}


\bibitem[Hil03]{Hi1} Hillman, J. A.: Tits alternatives and low
dimensional topology.  {\em J. Math. Soc. Japan\/} \textbf{55} (2003),
no. 2, 365--383.

\bibitem[Hi04]{Hi2} Hillman, J. A.: An indecomposable PD$_3$-complex:
II, {\em Algebraic and Geometric Topology\/} \textbf{4} (2004),
1103-1109.


\bibitem[IK04]{IK} Ivanov, S.; Katz, M.: Generalized degree and
optimal Loewner-type inequalities, {\em Israel J. Math.}  \textbf{141}
(2004), 221--233.  Available at the site
\texttt{arXiv:math.DG/0405019}


\bibitem[Ka83]{Ka1} Katz, M.: The filling radius of two-point
homogeneous spaces, {\em J. Diff. Geom.} \textbf{18} (1983), 505--511.


\bibitem[Ka88]{Ka1b} Katz, M.: The first diameter of~$3$-manifolds of
positive scalar curvature.  {\em Proc. Amer. Math. Soc.} \textbf{104}
(1988), no. 2, 591--595.


\bibitem[Ka06]{Ka4} Katz, M.: Systolic inequalities and Massey
products in simply-connected manifolds.  {\em Israel J. Math.} (2007).
See \texttt{arXiv:math.DG/0604012}


\bibitem[Ka07]{SGT} Katz, M.: Systolic geometry and topology.  With an
appendix by Jake P. Solomon.  {\em Mathematical Surveys and
Monographs}, \textbf{137}.  American Mathematical Society, Providence,
RI, 2007.



\bibitem[KL05]{KL} Katz, M.; Lescop, C.: Filling area conjecture,
optimal systolic inequalities, and the fiber class in abelian covers.
Geometry, spectral theory, groups, and dynamics, 181--200, {\em
Contemp. Math.} \textbf{387}, Amer. Math. Soc., Providence, RI, 2005.
See \texttt{arXiv:math.DG/0412011}


\bibitem[KR06]{KR1} Katz, M.; Rudyak, Y.: Lusternik-Schnirelmann
category and systolic category of low dimensional manifolds.  {\em
Communications on Pure and Applied Mathematics}, \textbf{59} (2006),
no.~10, 1433-1456.  Available at the site
\texttt{arXiv:math.DG/0410456}


\bibitem[KRS06]{KRS} Katz, M.; Rudyak, Y.: Sabourau, S.: Systoles of
2-complexes, Reeb graph, and Grushko decomposition.  {\em
International Math. Research Notices}, 2006 (2006).  Art.~ID 54936,
pp.~1--30.  See \texttt{arXiv:math.DG/0602009}




\bibitem[KS05]{KS2} Katz, M.; Sabourau, S.: Entropy of systolically
extremal surfaces and asymptotic bounds, {\em Ergo. Th. Dynam. Sys.},
\textbf{25} (2005), no.~4, 1209-1220.  See
\texttt{arXiv:math.DG/0410312}


\bibitem[KS06]{KS1} Katz, M.; Sabourau, S.: Hyperelliptic surfaces are
Loewner, {\em Proc. Amer. Math. Soc.} \textbf{134} (2006), no.~4,
1189-1195.  Available at the site \texttt{arXiv:math.DG/0407009}



\bibitem[KS06b]{KS3} Katz, M.; Sabourau, S.: An optimal systolic
inequality for CAT(0) metrics in genus two.  {\em Pacific J. Math.}
\textbf{227} (2006), no.~1, 95-107.  See
\texttt{arXiv:math.DG/0501017}


\bibitem[KSV06]{KSV} Katz, M.; Schaps, M.; Vishne, U.: Logarithmic
growth of systole of arithmetic Riemann surfaces along congruence
subgroups.  {\em J. Differential Geom.}  (2007).  Available at
\texttt{arXiv:math.DG/0505007}


\bibitem[LM89]{LM} Lawson, H. B.; Michelsohn, M.-L.: Spin geometry.
{\em Princeton Mathematical Series\/} \textbf{38}.  Princeton
University Press, Princeton, NJ, 1989.


\bibitem[LS34]{LS} Lusternik, L.; Schnirelmann, L.: M\'ethodes
topologiques dans les probl\`emes variationnels. Hermann, Paris 1934.

\bibitem[MH73]{MH} Milnor, J.; Husemoller, D.: Symmetric bilinear
forms.  Springer, 1973.


\bibitem[OR01]{OR} Oprea, J.; Rudyak, Y.: Detecting elements and
Lusternik-Schnirelmann category of 3-manifolds.
Lusternik-Schnirelmann category and related topics (South Hadley, MA,
2001), 181--191, {\em Contemp. Math.}, \textbf{316}, Amer. Math. Soc.,
Providence, RI, 2002.



\bibitem[We05]{We} Wenger, S.: Isoperimetric inequalities of
Euclidean type in metric spaces.  {\em Geom. Funct. Anal.} \textbf{15}
(2005), no.~2, 534--554.


\bibitem[Wh99]{Wh} White, B.: The deformation theorem for flat chains.
{\em Acta Math.} \textbf{183} (1999), no. 2, 255--271.


\bibitem[Wr74]{Wr} Wright, A. H.: Monotone mappings and degree one
mappings between~$PL$ manifolds.  Geometric topology (Proc. Conf.,
Park City, Utah, 1974), pp.~441--459.  {\em Lecture Notes in Math.}
\textbf{438}, Springer, Berlin, 1975.


\end{thebibliography}
\end{document}